\documentclass[10pt, titlepage]{amsart}

\usepackage{ae} 
\usepackage[T1]{fontenc}
\usepackage[cp1250]{inputenc}
\usepackage{amsmath}
\usepackage{amssymb, amsfonts,amscd,verbatim}
\usepackage[dvips]{graphicx}
\usepackage{latexsym}
\usepackage{indentfirst}
\usepackage{latexsym}

\usepackage{graphicx}
\usepackage{verbatim}   
\usepackage{color}      
\usepackage{subfigure}  

\usepackage{amsmath}    

\theoremstyle{plain}
\newtheorem{Prop}{Proposition}[section]
\newtheorem{Thm}[Prop]{Theorem}
\newtheorem{Cor}[Prop]{Corollary}

\theoremstyle{definition}
\newtheorem{Def}[Prop]{Definition}

\theoremstyle{remark}

\def\int{\mathop{\roman{int}}}

\def\1{^{-1}}

\def\NN{{\mathbf N}}

\def\dim{\text{dim}}

\def\diam{\text{diam}}
\def\dist{\text{dist}}

\def\UU{{\mathcal U}}

\def\asdim{\mathrm{asdim}}

\def\dim{\mathrm{dim}}
\def\diam{\mathrm{diam}}
\def\dokaz{{\bf Proof. }}
\def\edokaz{\hfill $\blacksquare$}

\numberwithin{equation}{section}

\input pstricks.tex
\input xy
\xyoption{all}

\begin{document}
\title[
Asymptotic dimension, Property A, and Lipschitz maps
]%
   {Asymptotic dimension, Property A, and Lipschitz maps}

\author{M.~Cencelj}
\address{IMFM,
Univerza v Ljubljani,
Jadranska ulica 19,
SI-1111 Ljubljana,
Slovenija }
\email{matija.cencelj@guest.arnes.si}

\author{J.~Dydak}
\address{University of Tennessee, Knoxville, TN 37996, USA}
\email{dydak@math.utk.edu}

\author{A.~Vavpeti\v c}
\address{Fakulteta za Matematiko in Fiziko,
Univerza v Ljubljani,
Jadranska ulica 19,
SI-1111 Ljubljana,
Slovenija }
\email{ales.vavpetic@fmf.uni-lj.si}

\date{ \today
}


\thanks{This research was supported by the Slovenian Research
Agency grants P1-0292-0101 and J1-2057-0101.}
\thanks{The second-named author was partially supported
by MEC, MTM2006-0825.}

\begin{abstract}

It is well-known that a paracompact space $X$ is of covering dimension $n$
if and only if any map $f\colon X\to K$ from $X$ to a simplicial complex $K$ can be pushed into its $n$-skeleton
$K^{(n)}$. We use the same idea to define dimension in the coarse category.
It turns out the analog of maps $f\colon X\to K$ is related to asymptotically Lipschitz maps,
the analog of paracompact spaces are spaces related to Yu's  Property A,
and the dimension coincides with Gromov's asymptotic dimension.

\end{abstract}

\maketitle

\medskip
\medskip
\tableofcontents
\section{Introduction}

It is well-known (see \cite{Dyd}) that the covering dimension $\dim(X)$
of a paracompact space can be defined as the smallest integer $n$
with the property that any commutative diagram
$$
\xymatrix{ A \ar_{i}[d] \ar^{g}[r]&   K^{{(n)}} \ar_{i}[d]  \\
X \ar^{f}[r]&   K 
}
$$
has a filler $h$
$$
\xymatrix{ A \ar_{i}[d] \ar^{g}[r]&   K^{{(n)}} \ar_{i}[d]  \\
X \ar^{h}[ur]\ar^{f}[r]&   K 
}
$$

Here $A$ is any closed subset of $X$, $K$ is any simplicial complex with the metric topology,
$K^{(n)}$ is the $n$-skeleton of $K$, and $i\colon A\to X$, $i\colon K^{(n)}\to K$ are inclusions.
By saying $h$ is a {\bf filler} we mean $h\vert A=g$ and, since we cannot insist on
$i\circ h=f$, we require $h(x)\in\Delta$ whenever $f(x)\in\Delta$ for any simplex $\Delta$
of $K$.

We want to prove the analog of the above result in the coarse category (see \cite{Roe lectures}).
The first issue is to find the analog of continuous maps $f\colon X\to K$.

As seen in \cite{Dyd} the optimal way to define paracompact spaces $X$
is as follows: for each open cover $\mathcal{U}$ of $X$ there is
a simplicial complex $K$ and a continuous map $f\colon X\to K$ such that
the family $\{f^{-1}(st(v))\}_{v\in K^{(0)}}$ refines $\mathcal{U}$.

In coarse geometry one replaces open covers by uniformly bounded covers,
refinement of covers by enlargement of covers, and continuous maps by either
bornologous functions or by $(\lambda,C)$-Lipschitz functions.
Since any function to a bounded space is bornologous, we will go with 
$(\lambda,C)$-Lipschitz functions.

Here is an analog of paracompact spaces in coarse geometry:

\begin{Def}\label{ls-paracompactDef}
$X$ is {\bf large scale paracompact} (ls-paracompact for short)
if for each uniformly bounded cover $\mathcal{U}$ of $X$ and for all $\lambda,C > 0$
there is a $(\lambda,C)$-Lipschitz function $f\colon X\to K$
such that $\mathcal{V}:=\{f^{-1}(st(v))\}_{v\in K^{(0)}}$ is uniformly bounded
and $\mathcal{U}$ refines $\mathcal{V}$.
\end{Def}

To simplify \ref{ls-paracompactDef} let us introduce the following concept:

\begin{Def}
Given $\delta > 0$ and a simplicial complex $K$, a function $f\colon X \to K$
is called a {\bf $\delta$-partition of unity} if it is $(\delta,\delta)$-Lipschitz,
$\mathcal{V}:=\{f^{-1}(st(v))\}_{v\in K^{(0)}}$ is uniformly bounded,
and the Lebesgue number of $\mathcal{V}$
is at least $\frac{1}{\delta}$.
\end{Def}

\begin{Prop}
A metric space $X$ is large scale paracompact if and only if for each $\delta > 0$
there is a $\delta$-partition of unity $f\colon X\to K$.
\end{Prop}
\dokaz Suppose for each $\delta > 0$ there is a $\delta$-partition of unity
$f_\delta\colon X\to K_\delta$. Given a uniformly bounded cover $\mathcal{U}$ of $X$ and $\lambda,C > 0$ choose $\delta > 0$ such that $\delta < \min(\lambda,C)$
and the diameter of each element of $\mathcal{U}$ is at most $\frac{1}{\delta}$.
Notice $f_\delta$ is $(\lambda,C)$-Lipschitz and $\mathcal{U}$
refines $\mathcal{V}:=\{f_\delta^{-1}(st(v))\}_{v\in K_\delta^{(0)}}$, i.e. $X$
is ls-paracompact.
\par If $X$ is paracompact and $\delta > 0$ put $\lambda=C=\delta$
and $\mathcal{U}$ is the cover of $X$ by $\frac{1}{\delta}$-balls.
The barycentric map is a $\delta$-partition of unity.
\edokaz

Our main results are as follows:

\begin{Thm}\label{TheoremA}
 If $X$ is of asymptotic dimension at most $n\ge 0$, then for any $\epsilon > 0$
 there is $\delta > 0$ such that any commutative diagram
 $$
\xymatrix{ A \ar_{i}[d] \ar^{g}[r]&   K^{{(n)}} \ar_{i}[d]  \\
X \ar^{f}[r]&   K 
}
$$
where $f$ is a $\delta$-partition of unity
has a filler $h$
$$
\xymatrix{ A \ar_{i}[d] \ar^{g}[r]&   K^{{(n)}} \ar_{i}[d]  \\
X \ar^{h}[ur]\ar^{f}[r]&   K 
}
$$
that is an $\epsilon$-partition of unity.
\end{Thm}

\begin{Thm}\label{TheoremB}
 Suppose $n\ge 0$ and for any $\epsilon > 0$
 there is $\delta > 0$ such that any commutative diagram
 $$
\xymatrix{ A \ar_{i}[d] \ar^{g}[r]&   K^{{(n)}} \ar_{i}[d]  \\
X \ar^{f}[r]&   K 
}
$$
where $f$ is a $\delta$-partition of unity
has a filler $h$
$$
\xymatrix{ A \ar_{i}[d] \ar^{g}[r]&   K^{{(n)}} \ar_{i}[d]  \\
X \ar^{h}[ur]\ar^{f}[r]&   K 
}
$$
that is an $\epsilon$-partition of unity.
If $X$ is large scale paracompact, then its asymptotic dimension is at most $n$.
\end{Thm}

\section{Preliminaries}

Given a set $S$ of vertices by $\Delta(S)$ we mean the full complex
over $S$: the set of functions $f\colon S\to [0,1]$ with finite support
such that $\sum\limits_{s\in S}f(s)=1$. $\Delta(S)$ is a subset of $l^{1}(S)$,
the space of all functions $f\colon S\to \mathbb{R}$ such that the $l^{1}$-norm $
\| f\|_1=\sum\limits_{s\in S}|f(s)|$ 
of $f$ is finite. $\Delta(S)$ inherits the resulting metric from $l^{1}(S)$.

By a simplicial complex $K$ we mean a subcomplex of $\Delta(S)$ for some set
$S$ ($S$ could be larger than the set of vertices $K^{(0)}$ of $K$).

Any function $g\colon X\to K$ from a space $X$ to a simplicial complex
$K$ can be viewed as a point-finite partition of unity $\{g_v\}_{v\in K^{(0)}}$,
where $g_v(x):=f(x)(v)$.

Given a vertex $v\in K^{(0)}$ by the {\bf star} $st(v)$ of $v$ in $K$
we mean all $f\in K$ such that $f(v) > 0$. Geometrically, it is the union of interiors
of all simplices of $K$ containing $v$.

Given a cover $\UU=\{U_s\}_{s\in S}$ of a metric space $(X,d)$
there is a natural family of functions $\{f_s\}_{s\in S}$ associated to $\UU$;
$f_s(x)$ defined as $\dist(x,X\setminus U_s)$. To simplify matters, by {\bf the local Lebesgue number} $L_\UU(x)$
of $\UU$ at $x$ we mean $$\sup\{f_s(x)\mid s\in S\}$$ and by the (global) {\bf Lebesgue number}
$L(\UU)$ of $\UU$ we mean $$\inf\{L_\UU(x)\mid x\in X\}.$$
We are interested in covers with a positive Lebesgue number. For such covers
{\bf the local multiplicity} $m_\UU(x)$ can be defined as
$|T(x)|$, where $T(x)=\{s\in S\mid f_s(x) > 0\}$
and the global multiplicity $m(\UU)$ can be defined as $$\sup\{m_\UU(x)\mid x\in X\}.$$
If the multiplicity $m(\UU)$ is finite, then $\UU$ has a natural partition of unity
$\{\phi_s\}_{s\in S}$ associated to it:
$$\phi_s(x)=\frac{f_s(x)}{\sum\limits_{t\in S} f_t(x)}.$$
That partition can be considered as a {\bf barycentric map}
$\phi:X\to\NN(\UU)$ from $X$ to the {\bf nerve} of $\UU$. We
consider that nerve with $l_1$-metric. Recall $\NN(\UU)$
is a simplicial complex with vertices belonging to $\UU$
and $\{U_1,\ldots,U_k\}$ is a simplex in $\NN(\UU)$
if and only if $\bigcap\limits_{i=1}^kU_i\ne\emptyset$.

\begin{Prop}\label{LipOfBary}
$\phi:X\to\NN(\UU)$ is $\frac{4m(\UU)^2}{L(\UU)}$-Lipschitz.
\end{Prop}
\dokaz
Since each $f_s$ is
$1$-Lipschitz, $\sum\limits_{t\in S} f_t(x)$ is
$2m(\UU)$-Lipschitz and each $\phi_s$ is
$\frac{2m(\UU)}{L(\UU)}$-Lipschitz (use the fact that
$\frac{u}{u+v}$ is
$\frac{\max(Lip(u),Lip(v))}{\inf(u+v)}$-Lipschitz). Therefore
$\phi:X\to\NN(\UU)$ is $\frac{4m(\UU)^2}{L(\UU)}$-Lipschitz. See
\cite{BD} and \cite{BS} for more details and better estimates of
Lipschitz constants.
\edokaz

If $f\colon X\to K$ is a function from a metric space to a simplicial complex $K$,
we can talk about its {\bf Lebesgue number} as that of the cover
$\{f^{-1}(st(v))\}_{v\in K^{(0)}}$. Notice that multiplicity of that cover
being at most $n+1$ implies $f(X) \subset K^{(n)}$.

\begin{Prop}\label{LebesgueOfAMap}
If $f\colon X\to K^{(n)}$ is $(\lambda,C)$-Lipschitz,
then the Lebesgue number of 
$\mathcal{V}:=\{f^{-1}(st(v))\}_{v\in K^{(0)}}$ is at least 
$R=\frac{1-(n+1)\cdot C}{(n+1)\cdot\lambda}$.
\end{Prop}
\dokaz
Given $x\in X$ find a vertex $v$ of $K$ so that $f_v(x)\ge \frac{1}{n+1}$.
If $d_X(x,y) < R$, then $|f_v(x)-f_v(y)|\leq \|f(x)-f(y)\|_1\leq \lambda\cdot d_X(x,y)+C
< \lambda\cdot R+C\leq \frac{1}{n+1}$.
Therefore $f_v(y) > 0$ and $y\in f^{-1}(st(v))$.
\edokaz

\begin{Def}[Definition 1.2.2 in \cite{Willett}]
A function $f\colon X\to Y$ of metric spaces
is said to have {\bf $(R,\epsilon)$ variation} if $d_X(x,y)\leq R$
implies $d_Y(f(x),f(y)) < \epsilon$.
\end{Def}

\begin{Prop}\label{RVariationAndLip}
If $f\colon X\to K$ has $(R,\epsilon)$ variation, $\epsilon \leq 2$, and $K$ is a metric simplicial complex,
then $f$ is $(\frac{2-\epsilon}{R},\epsilon)$-Lipschitz.
\end{Prop}
\dokaz If $d(x,y)\leq R$, then $\|f(x)-f(y)\|_1 < \epsilon\leq \frac{2-\epsilon}{R}\cdot 
d(x,y)+\epsilon$. Otherwise $\|f(x)-f(y)\|_1\leq 2=\frac{2-\epsilon}{R}\cdot R+\epsilon\leq
\frac{2-\epsilon}{R}\cdot d(x,y)+\epsilon$.
\edokaz

Let us show large scale paracompactness is a coarse invariant (see \cite{Roe lectures}
for basic concepts of coarse geometry).
\begin{Prop}
If $g\colon X\to Y$ is a coarse embedding and $Y$ is large scale paracompact,
then $X$ is large scale paracompact.
\end{Prop}
\dokaz Suppose $2 > \delta > 0$. Put $R=\frac{2-\delta}{\delta}$
and find $S > 0$ such that $d_X(x,y)\leq R$ implies $d_Y(g(x),g(y))\leq S$.
Find $2 > \epsilon > 0$ with the property that $\epsilon < \frac{\delta}{2(S+1)}$
and the cover of $X$ by $g^{-1}(B(y,\frac{1}{\epsilon}))$, $y\in Y$,
has Lebesgue number at least $\frac{1}{\delta}$.

Choose an $\epsilon$-partition of unity $f\colon Y\to K$ and put $h=f\circ g$.
If $d_X(x,y)\leq R$, then $d_Y(g(x),g(y))\leq S$
and $\|h(x)-h(y)\|_1\leq S\cdot \epsilon+\epsilon<\delta$,
so $h$ is $(\delta,\delta)$-Lipschitz by \ref{RVariationAndLip}.

For each $x\in X$ there is $y\in Y$ satisfying $B(x,\frac{1}{\delta})
\subset g^{-1}(B(y,\frac{1}{\epsilon}))$.
As $B(y,\frac{1}{\epsilon})\subset f^{-1}(st(v))$ for some vertex $v$ of $K$,
we have $B(x,\frac{1}{\delta})\subset h^{-1}(st(v))$.

Finally the family $\{h^{-1}(st(v))\}_{v\in K^{(0)}}$ is uniformly bounded
as $g$ is a coarse embedding.
\edokaz

By $B(A,R)$ we mean the union of all balls $B(x,R)$ of radius $R$,
where $x\in A$.

\begin{Prop}\label{ExtendingPUOverANbhd}
Suppose $f\colon X\to K$ has $(R,\epsilon)$ variation and $n\ge 0$.
If $A\subset f^{-1}(K^{(n)})$, then there is $r\colon B(A,R)\to K^{(n)}$
that has $(R,(8n+5)\cdot\epsilon)$ variation and is a filler of 
$$
\xymatrix{ A \ar_{i}[d] \ar^{f}[r]&   K^{{(n)}} \ar_{i}[d]  \\
B(A,R) \ar^{f}[r]&   K 
}
$$
\end{Prop}
\dokaz Given $x\in B(A,R)$ enumerate the set $V(x)$ of vertices of $K$ with the property $f_v(x) > 0$
as $v(0), v(1),\ldots$ so that $f_{v(i)}(x)\ge f_{v(i+1)}(x)$.
If the number of such vertices is at most $n+1$ (in particular if $x\in A$), we define $r(x)=f(x)$.
Otherwise we put $r_{v(0)}(x)=f_{v(0)}(x)+\sum\limits_{k=n+1}^\infty f_{v(k)}(x)$,
$r_{v(i)}(x)=f_{v(i)}(x)$ for $0 < i\leq n$, and $r_v(x)=0$ if $f_v(x)=0$.

Clearly, $r$ is a filler of $$
\xymatrix{ A \ar_{i}[d] \ar^{f}[r]&   K^{{(n)}} \ar_{i}[d]  \\
B(A,R) \ar^{f}[r]&   K 
}
$$
Also, if $x\in B(A,R)\setminus A$, then there is $x^\prime\in A$ with $d(x,x^\prime) < R$,
so $\|f(x)-f(x^\prime)\|_1 < \epsilon$. 

If $V(x^\prime)$ is a subset of $\widetilde V(x)=\{v(0),\ldots,v(n)\}$, then 
\begin{align*}
\sum\limits_{k=n+1}^\infty f_{v(k)}(x) &\le\sum\limits_{k=n+1}^\infty f_{v(k)}(x)+\sum\limits_{v\not\in V(x)-\widetilde V(x)} |f_{v}(x)-f_{v}(x^\prime)|=\\
&=\sum\limits_{v\in K^{(0)}} |f_{v(k)}(x)-f_{v(k)}(x^\prime)|=\\
&=\|f(x)-f(x^\prime)\|_1<\epsilon
\end{align*}
and 
\begin{align*}
\|f(x)-r(x)\|_1 &= \sum\limits_{k=0}^\infty |f_{v(k)}(x)-r_{v(k)}(x)|=\\
&=|f_{v(0)}(x)-r_{v(0)}(x)|+\sum\limits_{k=n+1}^\infty f_{v(k)}(x)=\\
&=2\sum\limits_{k=n+1}^\infty f_{v(k)}(x)<2\cdot\epsilon.
\end{align*}

Otherwise, there is $i\leq n$ so that $v(i)\notin V(x')$ implying
$f_{v(i)}(x)<\epsilon$ and therefore $f_{v(k)}(x) <\epsilon$ for all $k > i$.
For $k>i$
\begin{align*}
f_{v(k)}(x^\prime)&=f_{v(k)}(x^\prime)-f_{v(k)}(x)+f_{v(k)}(x)\le\\
&\le \|f(x^\prime)-f(x)\|_1+f_{v(k)}(x)\le 2\cdot\epsilon
\end{align*}
Thus 
\begin{align*}
\sum\limits_{k=n+1}^\infty f_{v(k)}(x) &< \sum\limits_{k=n+1}^\infty |f_{v(k)}(x)-f_{v(k)}(x^\prime)|+\sum\limits_{k=n+1}^\infty f_{v(k)}(x^\prime)<\\
(2n+1)\cdot \epsilon
\end{align*}
and 
\begin{align*}
\|f(x)-r(x)\|_1=2\sum\limits_{k=n+1}^\infty f_{v(k)}(x) < 2(2n+1)\cdot\epsilon.
\end{align*}

Finally, if $d(x,y) \leq R$, then $\|f(x)-r(x)\|_1< 2(2n+1)\cdot\epsilon$
and $\|f(y)-r(y)\|_1 < 2(2n+1)\cdot\epsilon$ resulting in $\|r(x)-r(y)\|_1 <  (8n+5)\cdot\epsilon$
as $\|f(x)-f(y)\|_1 <\epsilon$.
\edokaz

\section{Asymptotic dimension}

\begin{Prop}
Every metric space $X$ of finite asymptotic dimension
is large scale paracompact.
\end{Prop}
\dokaz
Let $\asdim(X)=n < \infty$ and $\delta > 0$.
Given $R > 0$ choose a uniformly bounded cover open $\UU_R=\{U_{s}\}_{s\in S}$
of multiplicity at most $n+1$ and Lebesgue number at least $R$.
We want the corresponding barycentric map
 $f\colon X\to \NN(\UU_R)^{(n)}$ to be a $\delta$-partition of unity.
Since $f^{-1}(st(s))=U_s$, we need $R > \frac{1}{\delta}$.
Also, by \ref{LipOfBary} $f$ is $\frac{(n+1)^2}{4R}$-Lipschitz,
so choosing $R$ at least $\frac{(n+1)^2}{4\cdot \delta}$ makes
$f$  a $\delta$-partition of unity.
\edokaz

\begin{Prop}\label{AsdimInTermsOfDeltaMaps}
 Suppose $X$ is a metric space and $n\ge 0$.  If for each $\delta > 0$
 there is a set $S$ and a $(\delta,\delta)$-Lipschitz map $f\colon X\to \Delta(S)^{(n)}$
 such that the family $\{f^{-1}(st(v))\}_{s\in S}$ is uniformly bounded, then $X$ has asymptotic dimension at most $n$.
\end{Prop}
\dokaz
If $f\colon X\to \Delta(S)^{(n)}$ is a $(\delta,\delta)$-Lipschitz map, then
$\{f^{-1}(st(s))\}_{s\in S}$ is of multiplicity at most $n+1$
and its Lebesgue number is at least $R=\frac{1-(n+1)\cdot \delta}{(n+1)\cdot\delta}$
by \ref{LebesgueOfAMap}. As $R$ can be made arbitrarily large, the asymptotic
dimension of $X$ is at most $n$.
\edokaz

{\bf Proof of Theorem \ref{TheoremA}:}
$h$ is going to be constructed as $h=\alpha\cdot r+(1-\alpha)\cdot \beta$,
where 
\begin{itemize}
\item[a.]
$r\colon B(A,R)\to K^{(n)}$ is a filler
of
$$
\xymatrix{ A \ar_{i}[d] \ar^{g}[r]&   K^{{(n)}} \ar_{i}[d]  \\
B(A,R) \ar^{f}[r]&   K 
}
$$
that has $(R,\mu)$ variation (see \ref{ExtendingPUOverANbhd})
for some $R,\mu > 0$ to be determined later.
\item[b.]
$\alpha, 1-\alpha$ is the barycentric partition of unity determined by sets
$B(A,R)$ and $B(C,R)$, where $C=X\setminus B(A,R)$.
\item[c.]
$\beta\colon X\to K^{(n)}$ is a barycentric partition of unity
of Lebesgue number at least $R$ with the property that
$\beta_v(x) > 0$ implies $f_v(x) > 0$ for all $x\in X$.
\end{itemize}
One should think of $h$ as a function from $X$ to $l^1(K^{(0)})$
in which case the formula $h=\alpha\cdot r+(1-\alpha)\cdot \beta$
makes sense provided $r$ is extended arbitrarily outside of $B(A,R)$
(as $\alpha=0$ outside $B(A,R)$ any extension will do).

The above conditions ensure that $h$ is a filler of 
$$
\xymatrix{ A \ar_{i}[d] \ar^{g}[r]&   K^{{(n)}} \ar_{i}[d]  \\
X \ar^{f}[r]&   K 
}
$$
Indeed, if $h_v(x) > 0$ for some $x\in X$ and some vertex $v$ of $K$,
then either $x\in B(A,R)$ and $r_v(x) > 0$ (in which case $f_v(x) > 0$)
or $x\in B(C,R)$ and $\beta_v(x) > 0$ (again, $f_v(x) > 0$ in this case).

By \ref{RVariationAndLip} $r$ is $(\frac{2-\mu}{R},\mu)$-Lipschitz,
$\beta$ is $\frac{4(n+1)^2}{R}$-Lipschitz by \ref{LipOfBary}, and $\alpha$
is $\frac{32}{R}$-Lipschitz by \ref{LipOfBary} (notice the Lebesgue number of the cover
$\{B(A,R),B(C,R)\}$ of $X$ is at least $\frac{R}{2}$).
Let us estimate Lipschitz constants of $h$:
on $B(A,R)$ one has $\|\alpha(x)\cdot r(x)-\alpha(y)\cdot r(y)\|_1\leq
\|\alpha(x)\cdot (r(x)-r(y))\|_1+\|(\alpha(x)-\alpha(y))\cdot r(y)\|_1\leq
\|r(x)-r(y)\|_1+\|\alpha(x)-\alpha(y)\|_1\leq (\frac{2-\mu}{R}+\frac{32}{R})\cdot d(x,y)+\mu
\leq \frac{34}{R}\cdot d(x,y)+\mu$. More generally, a product
of $\lambda$-Lipschitz function from $X$ to $[0,1]$ and a $(\mu,D)$-Lipschitz function
from $X$ to $l^1(S)$ is $(\lambda+\mu,D)$-Lipschitz.
If $x\in B(A,R)$ and $y\in X\setminus B(A,R)=C$,
then $\|\alpha(x)\cdot r(x)\|_1=\alpha(x)=\frac{\dist(x,C)}{\dist(x,C)+\dist(x,X\setminus B(C,R))}\leq \frac{d(x,y)}{R/2}$.
Thus $\alpha\cdot r$ is $(\frac{34}{R},\mu)$-Lipschitz when considered on the whole of $X$.
Similarly, $(1-\alpha)\cdot\beta$ is $\frac{4(n+3)^2}{R}$-Lipschitz
and $h$ is $(\frac{4(n+5)^2}{R},\mu)$-Lipschitz.

If we start with a $\delta$-partition of unity $f$ and $R > 0$,
then $f$ has $(R,R\cdot\delta+\delta)$ variation, so
we can find $r$ that has $(R,(8n+5)\cdot(R+1)\cdot \delta)$ variation by \ref{ExtendingPUOverANbhd}.
Thus we put $\mu= (8n+5)\cdot(R+1)\cdot \delta$.

In view of \ref{LebesgueOfAMap} we need
\begin{equation}\label{InequalityA}
\frac{1-(n+1)\cdot \mu}{(n+1)\cdot \frac{4(n+5)^2}{R}}\ge \frac{1}{\epsilon}
\end{equation}

to ensure the Lebesgue number of $h$ is at least $\frac{1}{\epsilon}$.
Also, we want 
\begin{equation}\label{InequalityB}
\mu<\epsilon
\end{equation}

 and 
\begin{equation}\label{InequalityC}
\frac{4(n+5)^2}{R} < \epsilon
\end{equation}

so that $h$ is $(\epsilon,\epsilon)$-Lipschitz.

To be able to construct $\beta$ we choose, for each $R > 0$,
 an open covering $\mathcal{U}^R$ of $X$ of multiplicity at most $n+1$
and Lebesgue number at least $R$ such that each element of $\mathcal{U}^R$ is of diameter at most $S(R) > R$. Thus $S(R)$ is a function of $R$.

First thing we need
is $S(R)< \frac{1}{\delta}$. Indeed, given 
a $\delta$-partition of unity $f\colon X\to K$ and $R$ satisfying $S(R)< \frac{1}{\delta}$, each element $U$ of $\mathcal{U}^R$
is assigned a unique vertex $v(U)$ of $K$ so that $U\subset f^{-1}(st(v(U)))$.
Now we can define sets $U_v$ as the union of all $U\in\mathcal{U}^R$
satisfying $v(U)=v$. That results in a covering of $X$ indexed by vertices of $K$ of Lebesgue number at least $R$,
of multiplicity at most $n+1$
(if $x\in \bigcap\limits_{i=1}^{n+2}U(v_i)$, then there are elements $U_i$, $1\leq i\leq n+2$,
of $\mathcal{U}^R$ containing $x$, hence $v_i=v_j$ for some $i\ne j$), and $U_v\subset f^{-1}(st(v))$.
The resulting barycentric partition of unity $\beta$ has Lebesgue number at least $R$.

If we consider $R=k$, $\mu= (8n+5)\cdot(R+1)\cdot \delta$, and $\delta=\frac{1}{k\cdot S(k)}$, it is clear that for sufficiently
large $k$ all inequalities \ref{InequalityA} - \ref{InequalityC} are satisfied.
\edokaz

{\bf Proof of Theorem \ref{TheoremB}:}

As $X$ is large scale paracompact, there is a $\delta$-partition of unity
$f\colon X\to \Delta(S)$ and $M > 0$ such that $\diam(f^{-1}(st(v))) < M$ for all $v\in S$.
Let $h\colon X\to \Delta(S)^{(n)}$ be an $\epsilon$-partition of unity and a push
of $f$. Apply \ref{AsdimInTermsOfDeltaMaps}.
\edokaz

\section{Property A}

In this section we investigate the relation of large scale paracompactness
to the Property A of Yu (see \cite{NY} or \cite{Yu00}).

We will use the following definition of Property A taken from \cite{Willett}
(beware we do not assume $X$ is of bounded geometry):

\begin{Def}
 $X$ has {\bf Property A} if for all $R,\epsilon > 0$ there is $M > 0$ and a partition of unity
 $\{\phi_{s}\}_{s\in S}$ with the following two properties:
 \begin{itemize}
\item [a.] $d(x,y)\leq R$ implies $\sum\limits_{s\in S}|\phi_{s}(x)-\phi_{s}(y)| < \epsilon$,
\item[b.] the diameter of the support of each $\phi_{s}$ is at most $M$.
\end{itemize}

\end{Def}

\begin{Prop}
 A metric space $X$ has Property A if and only if for each $\delta > 0$
 there is a set $S$ and a $(\delta,\delta)$-Lipschitz map $f\colon X\to \Delta(S)$
 such that the family $\{f^{-1}(st(v))\}_{s\in S}$ is uniformly bounded.
\end{Prop}
\dokaz
Suppose $X$ has Property A and $2 > \delta > 0$.
Put $R=\frac{2}{\delta}$ and choose $M > 0$ and a partition of unity
 $\{\phi_{s}\}_{s\in S}$ with the following two properties:
 \begin{itemize}
\item [a.] $d(x,y)\leq R$ implies $\sum\limits_{s\in S}|\phi_{s}(x)-\phi_{s}(y)| < \delta$,
\item[b.] the diameter of the support of each $\phi_{s}$ is at most $M$.
\end{itemize}
$\phi$ can be interpreted as a function $\phi\colon X\to \Delta(S)$
with $(R,\delta)$ variation.
By \ref{RVariationAndLip}, $\phi$ is $((2-\delta)/R,\delta)$-Lipschitz,
hence it is $(\delta,\delta)$-Lipschitz as $R=\frac{2}{\delta}$.
Notice $\phi^{-1}(st(v))$ equals the support of $\phi_{v}$.

Conversely, given $R,\epsilon > 0$ put $\delta=\frac{\epsilon}{R+1}$ and choose $\phi\colon X\to K$
that is $(\delta,\delta)$-Lipschitz and the diameter of $\phi^{-1}(st(v))$ is at most $M$
for some $M > 0$.
If $d(x,y)\leq R$, then $\|\phi(x)-\phi(y)\|_1\leq \delta\cdot R+\delta = \epsilon$.

If $\phi_{v}(x) > 0$ and $\phi_{v}(y) > 0$ then $x,y\in \phi^{-1}(st(v))$ and $d(x,y)\leq M$.
\edokaz

\begin{Cor}\label{lsParacompactHavePropertyA}
Each large scale paracompact space $X$ has Property A.
\end{Cor}

\begin{Thm}
If $(X,d)$ is a metric space of bounded geometry, then the following conditions are equivalent:
\begin{itemize}
\item[a.] $X$ has Property A,
\item[b.] $X$ is large scale paracompact.
\end{itemize}
\end{Thm}
\dokaz a)$\implies$b). As shown in \cite{Willett} for each $R,\epsilon > 0$ there is $S > 0$
and finite non-empty subsets $A_x\subset B(x,S)\times N$, $x\in X$,
such that $\frac{\vert  A_x\Delta A_y\vert}{\vert A_x\cap A_y\vert} < \epsilon$ if $d(x,y) \leq R$.
Given $1 > \delta > 0$ choose a natural number $M\ge 2$ such that every ball $B(x,\frac{1}{\delta})$
contains at most $M$ elements.
Pick $R > \frac{2-\delta}{\delta}+M$ and choose $S > M+\frac{1}{\delta}$
and finite non-empty subsets $A_x\subset B(x,S)\times N$, $x\in X$,
such that $\frac{\vert  A_x\Delta A_y\vert}{\vert A_x\cap A_y\vert} < \frac{\delta}{8M}$ if $d(x,y) \leq R$.

If $|A_x| \ge \frac{8M}{\delta}$ define $C_x=A_x\cup B(x,\frac{1}{\delta})\times \{1\}$,
otherwise define $C_x$ as $R(x)\times \{1\}$, where $R(x)$ is the $R$-component of $x$
(i.e., the set of all points $y$ that can be reached from $x$ by an $R$-chain).

If $d(x,y)\leq R$ and $|A_x|< \frac{8M}{\delta}$,
then $A_x=A_y$, as otherwise $\frac{\vert  A_x\Delta A_y\vert}{\vert A_x\cap A_y\vert}\ge
\frac{1}{|A_x|}> \frac{\delta}{8M}$. Thus the function $y\to A_y$ is constant over $R(x)$
and $R(x)$ must be contained in $B(x,2S)$ (otherwise $A_x\subset B(x,S)\times N$ and $A_y\subset B(y,S)\times N$
would have to be disjoint for some $y\in R(x)\setminus B(x,2S))$.
That proves $\frac{\vert C_x\Delta C_y\vert}{\vert C_x\cap C_y\vert }= 0$ 
and $C_x\subset R(x,2S)\times N$ in that case.

Suppose $d(x,y)\leq R$ and $|A_x| \ge  \frac{8M}{\delta}$.
Now $\frac{\vert C_x\Delta C_y\vert}{\vert C_x\vert }
\leq \frac{\vert A_x\Delta A_y\vert}{\vert A_x\vert }+\frac{2M}{\vert A_x\vert}
\leq \frac{\delta}{8M}+\frac{\delta}{4}< \frac{\delta}{3}$.

Define a partition of unity $f\colon X\to \Delta(X)$ by
$f_y(x)=\frac{\vert (\{x\}\times N)\cap C_y\vert}{\vert C_y\vert}$.
Notice $B(y,\frac{1}{\delta})\subset f^{-1}(st(y))\subset B(y,2S)$ for each $y\in X$. 
That means the Lebesgue number of $f$ is at least $\frac{1}{\delta}$
and the cover $\{f^{-1}(st(y))\}_{y\in X}$ is uniformly bounded.

Let us prove $f$ has $(R,\delta)$ variation which will show $f$ is a $\delta$-partition of unity
in view of \ref{RVariationAndLip}.

Notice $\| \vert C_x \vert\cdot f_x- \vert C_y \vert\cdot f_y\|_1\leq \vert C_x\Delta C_y\vert$
if $d(x,y)\leq R$. Without loss of generality assume $\vert C_x\vert \ge \vert C_y\vert$.
Therefore $0 \leq \frac{\vert C_x\vert}{\vert C_y\vert}-1\leq 
\frac{\vert C_x\Delta C_y\vert}{\vert C_y\vert}+\frac{\vert C_x\cap C_y\vert}{\vert C_y\vert}-1\leq \frac{\delta}{3}$.
Now $\| \frac{\vert C_x\Delta C_y\vert}{\vert C_y\vert}\cdot f_x-  f_y\|_1\leq 
\frac{\vert C_x\Delta C_y\vert}{\vert C_y\vert}\leq \frac{\delta}{3}$
and $\|f_x-  f_y\|_1\leq \| \frac{\vert C_x\vert}{\vert C_y\vert}\cdot f_x-f_x\|_1+
\| \frac{\vert C_x\vert}{\vert C_y\vert}\cdot f_x-  f_y\|_1\leq
\frac{\vert C_x\vert}{\vert C_y\vert}-1+\frac{\delta}{3} < \delta$.

\par b)$\implies$a) follows from \ref{lsParacompactHavePropertyA}.
\edokaz


\begin{thebibliography}{99}

\bibitem{BD}
G.Bell and A.Dranishnikov,
\emph{On asymptotic dimension of groups acting
on trees}, Geom. Dedicata {\bf 103} (2004), 89--101.

\bibitem{BS}
S.Buyalo and V.Schroeder, \emph{Hyperbolic dimension of metric
spaces}, arXiv:math.GT/0404525


\bibitem{Dyd}   J.Dydak, {\em    Partitions of unity},    Topology Proceedings
  27 (2003),  125--171.
 \par\noindent
 {\tt\verb+http://front.math.ucdavis.edu/math.GN/0210379+}

\bibitem{Grom}
M. Gromov, {\em Asymptotic invariants for infinite groups}, in
Geometric Group Theory, vol. 2, 1--295, G. Niblo and M. Roller,
eds., Cambridge University Press, 1993.

\bibitem{HR} 
N. Higson and J. Roe, {\em Amenable group actions and the Novikov conjecture},
J. Reine Agnew. Math. 519 (2000).
 

\bibitem{NY}
P. Nowak and G. Yu,
{\em What is ... Property A?},
Notices of the AMS Volume 55, Number 4, pp.474--475.

\bibitem{Roe lectures}
J. Roe, {\em Lectures on coarse geometry}, University Lecture
Series, 31. American Mathematical Society, Providence, RI, 2003.

 
\bibitem{Willett}
R. Willett, {\em Some notes on Property A}, arXiv:math/0612492v2 [math.OA]

\bibitem{Yu00}
G. Yu, {\em The coarse Baum-Connes conjecture for spaces which admit a uniform embedding into Hilbert space}, Inventiones 139 (2000), pp. 201--240.
\end{thebibliography}
\end{document}